\documentclass[10pt]{article}
\usepackage{mathrsfs}
\usepackage{amssymb}
\usepackage{amsmath}
\usepackage{graphicx}
\usepackage{geometry}
\usepackage{cite}
\usepackage{amsthm}
\newtheorem{theorem}{\scshape \mdseries  Theorem}[section]
\newtheorem{lemma}[theorem]{\scshape \mdseries  Lemma}
\newtheorem{coro}[theorem]{\scshape \mdseries  Corollary}

\newtheorem{defn}{Definition}[section]

\begin{document}

\title{\sf Graphs having extremal monotonic topological indices
with bounded vertex $k$-partiteness
\thanks{
      Supported by the key project of the Outstanding Young Talent Support Program of the University of Anhui Province (gxyqZD2016367), Project of Nature Science of Chizhou University(2016ZR008,2015ZRZ005,2015ZR005)
and the National Science Foundation of China under Grant No.11601006. }}
\author{$ Fang\ Gao~^{a},~~Duo\textrm{-}Duo \ Zhao~^{a} ,~~Xiao\textrm{-}Xin \ Li~^{a} ,~~Jia\textrm{-}Bao \ Liu~^{b,}
\thanks{Corresponding author.  E-mail address:
 gaofang@czu.edu.cn(F. Gao), liujiabaoad@163.com(J. Liu). }$\\
\\
{\small\it $^{a}$School of Mathematics and Computer Science,~ Chizhou University,~
 Chizhou 247000, China}\\
{\small  \it  $^{b}$ School of Mathematics and Physics,~ Anhui Jianzhu University,~
  Hefei 230601,  China}
   }
\date{}
\maketitle

\noindent {\bf Abstract}: The vertex $k$-partiteness $v_k(G)$ of graph $G$
is defined as the fewest number of vertices whose deletion from $G$ yields a $k$-partite graph.
In this paper, we introduce two concepts: monotonic decreasing topological index
and monotonic increasing topological index, and characterize the extremal graphs having
the minimum Wiener index, the maximum Harry index, the maximum reciprocal degree distance,
the minimum eccentricity distance sum, the minimum adjacent eccentric distance sum index,
the maximum connective eccentricity index, the maximum Zagreb indices among graphs
with a fixed number $n$ of vertices and fixed vertex $k$-partiteness, respectively.

%
%
%\noindent


%\vskip 0.1cm

\vspace{-0.2cm}

\section{Introduction}
\vspace{-0.2cm}

~~~~Let $G$ be a simple graph with vertex-set $V(G)$ and edge-set $E(G).$
An edge with end vertices $u$ and $v$ is denoted by $uv$ and we say that $u$ and $v$ are
adjacent or neighbors. The degree of a vertex $u \in V(G),$ denoted by $d_G (u),$ is the number of neighbors of $u.$
The distance between two vertices $u$ and $v,$ denoted by $d_G (u,v),$
is the length of a shortest path connecting them in $G,$ and the sum of distance related to $u$ is defined as
$D_G(u)=\sum_{v\in V(G)}d_G (u,v).$
The eccentricity of a vertex $u$ , denoted by $\varepsilon_G (u),$ is the maximum distance from $u$ to
any other vertex. In the paper, We omit the subscript $G$ from the notation.
The complement of $G,$ denoted by $\overline G,$
is the graph with $V(\overline G)=V(G)$ and $E(\overline G)=\{uv:uv\notin E(G)\}.$
A subgraph $H$ of $G$ is an induced subgraph if two vertices of $V(H)$ are adjacent
in $H$ if and only if they are adjacent in $G.$ Hence, an induced subgraph is determined by its vertex set.
The induced subgraph with vertex set $S$ is denoted by $\langle S \rangle.$
A bipartite graph is a graph $G$ whose vertex set composed of two disjoint set $X,Y$ such that
each edge has an end vertex in $X$ and the other one in $Y.$
A complete bipartite graph $K_{s,t}$ is a bipartite graph with $|X|=s,|Y|=t$ and such that
any two vertices $u\in X$ and $v \in Y$ are adjacent.
The complete graph of order $n$ is denoted by $K_n.$
The join of two vertex disjoint graphs $G_1$ and $G_2,$ denoted by $G_1 \vee G_2,$
is the graph obtained from the disjoint union $G_1\cup G_2$ by adding edges between each vertex of $G_1$ and each of $G_2.$
For terminology and notation not defined here we refer to \cite{Bondy}.

Molecular graphs are models of molecules in which atoms are represented by
vertices and chemical bonds are represented by edges of a graph.
Chemical graph theory is a branch of mathematical chemistry concerning the study of chemical graph.
A topological index (also known as molecular descriptor or graph invariant)
is a single number that can be used to characterize some property of molecule graphs.
Topological index is a graph theoretic property which is preserved by isomorphism.
The chemical information derived through topological index has been found useful in chemical documentation,
isomer discrimination, structure property correlations, etc \cite{ash}. A large number of graph invariants
of molecular graphs are studied in chemical graph theory, based on vertex degree or vertex distance\cite{gut2013,Ni}.

One of the oldest and well-known graph invariants is Wiener index, denoted by $W(G),$ which is introduced
by Wiener \cite{wie} in 1947 and defined as the sum of distance over all unordered vertex pairs in $G,$ i.e.
$$W(G)=\sum\limits_{\{u,v\}\subseteq V(G)} d(u,v).$$

Another graph invariant, defined in a fully analogous manner to Wiener index, is Harary index,
denoted by $H(G),$ which is defined as the sum of reciprocals of distances between all pairs of vertices in $G,$ i.e.
$$H(G)=\sum\limits_{\{u,v\}\subseteq V(G)} \frac{1}{d(u,v)}.$$

In 2012, Hua and Zhang \cite{hua2} proposed a modification of Harary index, reciprocal degree distance,
which can be seen as the weighted degree-sum version of Harary index, and is defined as
$$RDD(G)=\sum\limits_{\{u,v\}\subseteq V(G)}\frac{d(u)+d(v)}{d(u,v)}.$$

In order to analyze the  structure-dependency of total $\pi$-electron energy on the molecular,
Gutman and Trinajst\'{c} \cite{gut1975,gut1972} introduced the
Zagreb indices, and defined as
$$M_1(G)=\sum\limits_{uv\in E(G)} (d(u)+d(v))=\sum\limits_{u\in V(G)} d(u)^2,M_2(G)=\sum\limits_{uv\in E(G)} d(u)d(v).$$

Todeschini and Consonni \cite{tod2010} considered the multiplicative versions of
the Zagreb indices in 2010, which was defined as
$$\Pi_1(G)=\prod \limits_{u\in V(G)}d(u)^2,\Pi_2(G)=\prod \limits_{uv\in E(G)}d(u)d(v)=\prod \limits_{u\in V(G)}d(u)^{d(u)}.$$

$M_1(G)$ is called the first Zagreb index and $M_2(G)$ the second Zagreb index.
$\Pi_1(G)$ is called the first multiplicative Zagreb index and $\Pi_2(G)$ the second multiplicative Zagreb index.

Recently, the distance-based invariants including eccentricity have attracted more and more attention.
For example, the eccentricity distance sum \cite{gup2002} of a graph $G,$ is defined as
$$\xi^d(G)=\sum\limits_{uv\in E(G)}(\varepsilon(u)+\varepsilon(v))d(u,v)=\sum\limits_{u\in V(G)}\varepsilon(u)D(u).$$

The connective eccentricity index\cite{gup2000} of a graph $G,$ is defined as
$$\xi^{ce}(G)=\sum\limits_{uv\in E(G)}(\frac{1}{\varepsilon(u)}+\frac{1}{\varepsilon(v)})=\sum\limits_{u\in V(G)} \frac{d(u)}{\varepsilon(u)}.$$

The adjacent eccentric distance sum index \cite{sar2002} of a graph $G,$ is defined as
$$\xi^{ad}(G)=\sum\limits_{u\in V(G)} \frac{\varepsilon(u)D(u)}{d(u)}.$$

These indices all above have attracted extensive attention due to their wide applications
in physics, chemistry, graph theory, etc.

In order to measure how close a graph $G$ is to being a bipartite graph, S. Fallat and
Yi-Zheng Fan \cite{Fallat} introduced a parameter:
the vertex bipartiteness, which is denoted by $v_b(G)$ and defined as
the fewest number of vertices whose deletion yields a bipartite graph.
Let $m$ be a natural number such that $m \leq n-2.$  Let
$$\mathscr {G}_{n,m}= \Big\{G = (V(G),E(G)) : |V(G)| =n, v_b(G)\leq m \Big\}.$$

In 2016, I. Gutman et al.\cite{gut2016} determined the maximum Laplacian Estrada index
and the maximum signless Laplacian Estrada index in $\mathscr{G}_{n,m}.$
M. A. A. de Abreu et al. \cite{Freitas2016} determined graphs in $\mathscr {G}_{n,m}$
having maximum Laplacian-energy-like invariant and the incidence energy.
M. Robbiano et al.\cite{rob} identified the graphs in $\mathscr{G}_{n,m}$
with maximum spectral radius and maximum signless Laplacian spectral radius.
J.B. Liu and X.F. Pan \cite{liu} characterized the graph in $\mathscr {G}_{n,m}$
having the minimum Kirchhoff index. H. Li et al.\cite{hong2017} characterized some
the maximal connective eccentricity indices in $\mathscr{G}_{n,m}.$ Motivated from these results, it is natural and
interesting to study the extremal indices with bounded vertex $k$-partiteness.

A $k$-partite graph is a graph $G$ whose vertex set composed of $k$-disjoint
set $U_1,U_2,\ldots,U_k$ such that each edge has an end vertex in $U_i$
and the other one in $U_j,$ for any $1\leq i\neq j\leq k.$
The vertex $k$-partiteness of graph $G$ is defined as the fewest number of vertices
whose deletion from $G$ yields a $k$-partite graph, denoted by $v_k(G).$  Let
$$\mathscr {G}_{n,m,k}= \Big\{G = (V(G),E(G)) : |V(G)| =n, v_k(G)\leq m,m\leq n-k \Big\}.$$

In this paper, we introduce two concepts: monotonic decreasing topological index and
monotonic increasing topological index, and characterize the graphs with extremal monotonic topological
indices in $\mathscr {G}_{n,m,k}.$ As applications, we identify
the graphs in $\mathscr {G}_{n,m,k}$ with the minimum Wiener index, the maximum Harary index,
the maximum reciprocal degree distance,
the minimum eccentricity distance sum,
the minimum adjacent eccentric distance sum index,
the maximum connective eccentricity index, and the maximum Zagreb indices respectively.

\vspace{-0.3cm}
\section{Preliminaries}\vspace{-0.3cm}

~~~~For a subset $W\subset V(G),$ let $G-W$ be the subgraph of $G$ obtained by deleting
the vertices of $W$ together with the edges incident with them.
For a subset $E_1\subset E(\overline{G}),$
let $G+E_1$ be the graph obtained from $G$ by adding the edges of $E_1.$
If $E_1=\{e\},$ we denote by $G+e$ for simplicity.

\begin{defn}
For some topological index $TI(G)$ of a graph $G,$
 if $TI(G+e)<TI(G),$ we call it monotonic decreasing topological
index, denoted by $TI_-(G).$
\end{defn}

\begin{defn}
For some topological index $TI(G)$ of a graph $G,$
if $TI(G+e)>TI(G),$ we call it monotonic increasing topological
index, denoted by $TI_+(G).$
\end{defn}

Let $e=uv$ be an edge of $\overline{G}.$
Adding $e=uv$ to $G$ does not increase distances,
while it does decrease at least one distance,
the distance between $u$ and $v$ is at least 2 in $G$ and 1 in $G+e,$ and the distance sum $D(u)$ of
vertex $u$ will decrease. At the same time, the adding of $e$ does not decrease vertices degree,
while it does increase the degree of $u$ and $v,$ and the eccentricity $\varepsilon(u)$ will not increase.

By the analysis above, we have the following lemmas immediately for a connected graph.
\begin{lemma}
 Let $G$ be a graph with $u,v \in V(G).$ If $uv \notin E(G),$ then
$W(G) > W(G+uv).$
\end{lemma}

\begin{lemma}
Let $G$ be a graph with $u,v \in V(G).$ If $uv \notin E(G),$ then $H(G) < H(G+uv).$
\end{lemma}

\begin{lemma}
Let $G$ be a graph with $u,v \in V(G).$ If $uv \notin E(G),$ then $RDD(G) < RDD(G+uv).$
\end{lemma}

\begin{lemma}
Let $G$ be a graph with $u,v \in V(G).$ If $uv \notin E(G),$ then $\xi^d(G) > \xi^d(G+uv).$
\end{lemma}

\begin{lemma}
Let $G$ be a graph with $u,v \in V(G).$ If $uv \notin E(G),$ then $\xi^{ce}(G) < \xi^{ce}(G+uv).$
\end{lemma}

\begin{lemma}
Let $G$ be a graph with $u,v \in V(G).$ If $uv \notin E(G),$ then $\xi^{ad}(G) > \xi^{ad}(G+uv).$
\end{lemma}

\begin{lemma}
Let $G$ be a graph with $u,v \in V(G).$ If $uv \notin E(G),$  then $M_1(G) < M_1(G+uv);\\
M_2(G) < M_2(G+uv); \Pi_1(G)<\Pi_1(G+uv); \Pi_2(G)<\Pi_2(G+uv).$
\end{lemma}

Hence, the Wiener index, the eccentricity distance sum index, the adjacent eccentric distance sum index are
 monotonic decreasing topological indices,
while the Harary index, the reciprocal degree distance, the connective eccentricity index, the
Zagreb indices are monotonic increasing topological indices.

\section{Main results}

~~~~In this section, we will characterize the graphs with extremal monotonic decreasing
topological index or monotonic increasing topological index in $\mathscr {G}_{n,m,k}.$

\begin{theorem}
Let $1 \leq m\leq n-k.$ Then there exist $k$ non-negative integers $s_1,s_2,\ldots,s_k$
satisfying $s_1+s_2+\ldots +s_k=n-m,$ such that $TI_-(G)\geq TI_-(\widehat{G})$ holds for
all graphs $G\in \mathscr {G}_{n,m,k},$ where
 $$  \widehat{G} =K_{m}\vee\big(\overline{K_{s_1}}\vee\overline{K_{s_2}}\vee\ldots\vee\overline{K_{s_k}}\big)\in \mathscr {G}_{n,m,k},$$
with equality holds if and only if $ G \cong \widehat{G}.$
\end{theorem}
{\bf Proof.}
Denote by $G^+$ the graph obtained from $G$ by adding an edge $uv,$ where $u,v\in V(G),uv\notin E(G).$
Let $\widehat{G}\in \mathscr {G}_{n,m,k}$  such that $TI_-(G)\geq TI_-(\widehat{G})$
for all $G\in \mathscr {G}_{n,m,k}.$
For simplicity, we denote by $k^*$ the $k$-partiteness of graph $\widehat G,$
then $k^*\leq m.$ Hence, there exist $k^*$ vertices $i_1,i_2,\ldots,i_{k^*} \in V(\widehat{G})$
such that$\widehat{G}-\{i_1,i_2,\ldots,i_{k^*}\}$ is a $k$-partite graph with $k$-partition $\{U_1,U_2,\ldots,U_k\}.$
Let $s_i = |U_i|,~1\leq i\leq k,$ hence $n = s_1 + \ldots +s_k +{k^*}.$

Firstly, we claim that
$$\widehat{G}-\{i_1,i_2,\ldots ,i_{k^*}\}=K_{s_1 ,s_2, \ldots ,s_k}
=\overline{K_{s_1}}\vee \overline{K_{s_2}}\vee\ldots\vee \overline{K_{s_k}}.$$
Otherwise, there exist two vertices $u \in U_{s_i}$ and $v \in U_{s_j} ,1\leq i\neq j\leq k,$
which are not adjacent in $\widehat{G}.$
Denote by $\widehat{G}^+$ the graph $\widehat{G} + uv,$ obviously, $\widehat{G}^+\in \mathscr {G}_{n,m,k}.$
By the definition of monotonic decreasing topological index, we get
$$TI_-(\widehat{G}^+)< TI_-(\widehat{G}),$$
which is a contradiction.

Secondly, we claim that$\langle \{i_1,i_2,\ldots ,i_{k^*}\} \rangle=K_{k^*}.$
Otherwise, there exist two vertices $u,v$ are not adjacent,
where $u, v \in \{i_1,i_2,\ldots ,i_{k^*}\}.$
By connecting the vertices $u$ and $v,$
we arrive at a new graph $\widehat{G} + uv,$ obviously, $\widehat{G} + uv \in\mathscr {G}_{n,m,k}.$
By the definition of monotonic decreasing topological index,
we get $$TI_-(\widehat{G}+uv)< TI_-(\widehat{G}),$$
a contradiction again.

Thirdly, using a similar method, we can get
$$\widehat{G} =K_{k^*}\vee\big(\overline{K_{s_1}}\vee\overline{K_{s_2}}\vee\ldots\vee\overline{K_{s_k}}\big).$$

Finally, we prove $k^*= m.$ If $k^*\leq m-1,$ then $s_1+s_2+\ldots +s_k=n-k^* \geq n-m+1>n-m\geq k,$
thus $s_1+s_2+\ldots +s_k>k.$
Without loss of generality, we assume that $s_1\geq 2.$
Note that the graph $\widehat{G}=K_{k^*}\vee\big(\overline{K_{s_1}}\vee\overline{K_{s_2}}\vee\ldots\vee\overline{K_{s_k}}\big)$
contains the subgraph $K_{k^*+1},$
which is an induced subgraph of $\widehat{G}$ with vertex set $S=\{u,i_1,i_2,\ldots,i_{k^*}\},$
where $u$ is a vertex of $K_{s_1}.$
By connecting $u$ and other vertices of $\big(\overline{K_{s_2}}\vee\ldots\vee\overline{K_{s_k}}),$
we get a new graph $\widetilde{G} =K_{k^*+1}\vee (\overline{K_{s_1-1}}\vee\overline{K_{s_2}}
\vee\ldots\vee\overline{K_{s_k}}) \in \mathscr {G}_{n,m,k},$
which has $s_1-1\geq 1$ edges more than the graph $\widehat{G}.$
By the definition of monotonic decreasing topological index, we get
$$TI_-(\widetilde{G})> TI_-(\widehat{G}),$$
which is obviously another contradiction.

Therefore
$$ \widehat{G} =K_m\vee\big(\overline{K_{s_1}}\vee\overline{K_{s_2}}\vee\ldots\vee\overline{K_{s_k}}\big).$$
The proof of the theorem is complete.\qed

For a monotonic increasing topological index, we have a similar result.
\begin{theorem}
Let  $1 \leq m\leq n-k.$ Then there exist $k$ non-negative integers $s_1,s_2,\ldots,s_k$
satisfying $s_1+s_2+\ldots+s_k=n-m,$ such that $TI_+(G)\leq TI_+(\widehat{G})$
holds for all graphs $G\in \mathscr {G}_{n,m,k},$  where
$$\widehat{G} =K_{m}\vee\big(\overline{K_{s_1}}\vee\overline{K_{s_2}}\vee\ldots\vee
\overline{K_{s_k}}\big)\in \mathscr {G}_{n,m,k},$$
with equality holds if and only if $ G \cong \widehat{G}.$
\end{theorem}

\section{Applications}

~~~~As some applications of Theorem $3.1$ and $3.2,$ we will investigate
the minimum Wiener index, the maximum Harary index, the maximum
reciprocal degree distance, the minimum eccentricity distance sum index,
the minimum adjacent eccentric distance sum index, the maximum connective eccentricity index, the
maximum Zagreb indices in $\mathscr{G}_{n,m,k},$ respectively.
Assume that $n-m=sk+t,$ where $s$ and $t$ are non-negative integers such that $0\leq t < k.$

\subsection{The minimum Wiener index in $\mathscr{G}_{n,m,k}$}

\begin{theorem}
Let $G$ be a connected simple graph of order $n$ with vertex $k$-partiteness $v_k(G) \leq m,$
where $1\leq m \leq n-k.$ Then
  $$ W(G)\geq\frac{n^2-m}{2}+\frac{(n-m)(s-2)}{2}+\frac{t(s+1)}{2} $$
with equality holds if and only if $G\cong K_m\vee((k-t)\overline {K_s}\vee t\overline {K_{s+1}}),$\\
where $(k-t)\overline {K_s}=\underbrace{\overline {K_s}\vee\ldots\vee\overline {K_s}}_{k-t},$ and
$t\overline {K_{s+1}}=\underbrace{\overline {K_{s+1}}\vee\ldots\vee\overline {K_{s+1}}}_t.$
\end{theorem}
{\bf Proof.} Since the Wiener index is a monotonic decreasing
topological index, by Theorem $3.1,$ there exist $k$ non-negative integers $s_1,s_2,\ldots,s_k$
satisfying $s_1+s_2+\ldots +s_k=n-m,$
such that $\widehat{G}=K_m\vee\big(\overline{K_{s_1}}\vee\overline{K_{s_2}}\vee\ldots\vee\overline{K_{s_k}}\big)$
is the graph in $\mathscr {G}_{n,m,k}$ with the minimum Wiener index.

In the following, we will determine the values of $s_1,s_2,\ldots,s_k.$\\
$W(K_m\vee(\overline{K_{s_1}} \vee\overline{K_{s_2}} \vee\ldots\vee\overline{K_{s_k}}))$\\
=$\sum\limits_{\{u,v\}\subseteq V(\widehat{G})}d(u,v)$\\
=$C_m^2\times 1 +\sum\limits_{i=1}^{k} C_{s_i}^2 \times 2+
\sum\limits_{i=1}^{k} m\times s_i \times 1+\sum\limits_{1\leq i<j\leq k}s_is_j \times 1$\\
=$\sum\limits_{i=1}^{k}s_i(s_i-1)+\sum\limits_{1\leq i<j\leq k}s_is_j+C_m^2+m(n-m).$

We claim that $\widehat{G}$ is the graph in $\mathscr {G}_{n,m,k}$ with the minimum Wiener index
when $s_i = s_j $ or $|s_i-s_j|=1,$ for all $1\leq i,j\leq k.$
If there exist $1\leq i,j\leq k,$ such that $|s_i-s_j|\geq 2.$
Wihtout loss of generality, we assume that $s_1-s_2\geq2.$
By moving one vertex from the part of $\overline{K_{s_1}}$ to the part of  $\overline{K_{s_2}},$  we get a new graph
$\widetilde{G}=K_m\vee(\overline {K_{s_1-1}} \vee\overline{K_{s_2+1}}\vee\ldots\vee\overline{K_{s_k}}),$
which is also in $\mathscr{G}_{n,m,k}.$ Then\\
$W(\widehat{G})-W(\widetilde{G})\\
=s_1(s_1-1)+s_2(s_2-1)+\sum\limits_{i=2}^{k}s_1s_i+\sum\limits_{i=3}^{k}s_2s_i\\
  -(s_1-1)(s_1-2)-(s_2+1)s_2-(s_1-1)(s_2+1)-\sum\limits_{i=3}^{k}(s_1-1)s_i-\sum\limits_{i=3}^{k}(s_2+1)s_i\\
=s_1-s_2-1.$

Since $s_1-s_2\geq2,$ then $W(\widehat{G}) > W(\widetilde{G}),$
which is a contradiction.

Therefore, for any $1\leq i,j\leq k,$ then $ |s_i-s_j|\leq 1.$ For $s_1+s_2+\ldots +s_k=n-m=sk+t,$
we have $\overline{K_{s_1}} \vee\overline{K_{s_2}} \vee\ldots\vee\overline{K_{s_k}}=(k-t)\overline{K_s}\vee t\overline{K_{s+1}},$\\
where $(k-t)\overline{K_s}=\underbrace{\overline{K_s}\vee\ldots\vee\overline{K_s}}_{k-t},$
and $t\overline{K_{s+1}}=\underbrace{\overline{K_{s+1}}\vee\ldots\vee\overline{K_{s+1}}}_t.$

Then we obtain that\\
$W(\widehat{G})=\frac{n^2-m}{2}+\frac{(n-m)(s-2)}{2}+\frac{t(s+1)}{2}.$ \qed

\begin{coro}
Let $G\in\mathscr{G}_{n,m,2}.$  Then the following holds

(1) If $n-m$ is even, then
$$W(G) \geq\frac{3n^2+m^2-2mn-4n+2m}{4},$$

with equality holds if and only if
$$G\cong K_{m}\vee\big(\overline{K_{\frac{n-m}{2}}}\vee\overline{K_{\frac{n-m}{2}}}\big).$$

(2) If $n-m$ is odd, then
$$W(G) \geq\frac{3n^2+m^2-2mn-4n+2m+1}{4},$$

with equality holds if and only if
$$G\cong K_{m}\vee\big(\overline{K_{\frac{n-m+1}{2}}}\vee\overline{K_{\frac{n-m-1}{2}}}\big).$$

\end{coro}

\subsection{The maximum Harry index in $\mathscr{G}_{n,m,k}$}

\begin{theorem}

Let $G$ be a connected simple graph of order $n$ with vertex $k$-partiteness $v_k(G) \leq m,$
where $1\leq m \leq n-k.$ Then
  $$H(G)\leq\frac{n^2-m}{2}-\frac{(n-m)(s+1)}{4}-\frac{t(s+1)}{4}$$
with equality holds if and only if $G\cong K_m\vee((k-t)\overline{K_s}\vee t\overline{K_{s+1}}),$\\
where $(k-t)\overline{K_s}=\underbrace{\overline{K_s}\vee\ldots\vee\overline{K_s}}_{k-t},$
and $t\overline{K_{s+1}}=\underbrace{\overline{K_{s+1}}\vee\ldots\vee\overline{K_{s+1}}}_t.$
\end{theorem}

 {\bf Proof.} Since the Harary index is a monotonic increasing
topological index, by Theorem $3.2,$ there exist $k$ non-negative integers $s_1,s_2,\ldots,s_k$
satisfying $s_1+s_2+\ldots+s_k=n-m,$
such that $\widehat{G}=K_m\vee\big(\overline{K_{s_1}}\vee\overline{K_{s_2}}\vee\ldots\vee\overline{K_{s_k}}\big)$
is the graph in $\mathscr {G}_{n,m,k}$ with the maximum Harary index.

In the following, we will determine the values of $s_1,s_2,\ldots,s_k.$

$H(K_m\vee(\overline{K_{s_1}} \vee\overline{K_{s_2}} \vee\ldots\vee\overline{K_{s_k}}))$\\
=$\sum\limits_{\{u,v\}\subseteq V(\widehat{G})}\frac{1}{d(u,v)}$\\
=$ C_m^2\times1+\sum\limits_{i=1}^{k} C_{s_i}^2\times \frac{1}{2}
+\sum\limits_{i=1}^{k} m\times s_i\times 1+\sum\limits_{1\leq i<j\leq k}s_is_j\times 1$\\
=$\sum\limits_{i=1}^{k}\frac{s_i(s_i-1)}{4}+\sum\limits_{1\leq i<j\leq k}s_is_j+C_m^2+m(n-m).$

We claim that $\widehat{G}$ is the graph in $\mathscr {G}_{n,m,k}$ with the maximum Harary index
when $s_i = s_j $ or $|s_i-s_j|=1,$ for all $1\leq i,j\leq k.$
If there exist $1\leq i,j\leq k,$ such that $|s_i-s_j|\geq 2.$
Wihtout loss of generality, we assume that $s_1-s_2\geq2.$
By moving one vertex from the part of $\overline{K_{s_1}}$ to the part of  $\overline{K_{s_2}},$
we get a new graph
$\widetilde{G}=K_m\vee(\overline{K_{s_1-1}} \vee\overline{K_{s_2+1}}\vee\ldots\vee\overline{K_{s_k}}),$
which is also in $\mathscr{G}_{n,m,k}.$

Then we obtain that\\
$H(\widehat{G})-H(\widetilde{G})$\\
=$\frac{s_1(s_1-1)}{4}+\frac{s_2(s_2-1)}{4}+\sum\limits_{i=2}^{k}s_1s_i+\sum\limits_{i=3}^{k}s_2s_i\\
-\frac{(s_1-1)(s_1-2)}{4}-\frac{(s_2+1)s_2}{4}-(s_1-1)(s_2+1)\\
-\sum\limits_{i=3}^{k}(s_1-1)s_i-\sum\limits_{i=3}^{k}(s_2+1)s_i$\\
=$-\frac{1}{2}(s_1-s_2-1).$

Since $s_1-s_2\geq2,$ then $H(\widehat{G})< H(\widetilde{G}),$
which is a contradiction.

Therefore, for any $1\leq i,j\leq k,$ then $ |s_i-s_j|\leq 1.$ For $s_1+s_2+\ldots +s_k=n-m=sk+t,$
we have $\overline{K_{s_1}} \vee\overline{K_{s_2}} \vee\ldots\vee\overline{K_{s_k}}=(k-t)\overline{K_s}\vee t\overline{K_{s+1}},$\\
where $(k-t)\overline{K_s}=\underbrace{\overline{K_s}\vee\ldots\vee\overline{K_s}}_{k-t},$ and
$t\overline{K_{s+1}}=\underbrace{\overline{K_{s+1}}\vee\ldots\vee\overline{K_{s+1}}}_t.$

Then we obtain that\\
$H(\widehat{G})=\frac{n^2-m}{2}-\frac{(n-m)(s+1)}{4}-\frac{t(s+1)}{4}.$\qed

\begin{coro}
Let $G\in\mathscr{G}_{n,m,2}.$ Then the following holds.

(1) If $n-m$ is even, then
$$H(G) \leq \frac{3n^2-m^2+2mn-2m-2n}{8},$$
with equality holds if and only if
$$G\cong K_{m}\vee\big(\overline{K_{\frac{n-m}{2}}}\vee\overline{K_{\frac{n-m}{2}}}\big).$$

(2) If $n-m$ is odd, then
$$H(G) \leq \frac{3n^2-m^2+2mn-2m-2n-1}{8},$$
 with equality holds if and only if
 $$G\cong K_{m}\vee\big(\overline{K_{\frac{n-m+1}{2}}}\vee\overline{K_{\frac{n-m-1}{2}}}\big).$$

\end{coro}

\subsection{The maximum reciprocal degree distance in $\mathscr{G}_{n,m,k}$}

\begin{theorem}
Let $G$ be a connected simple graph of order $n$ with vertex $k$-partiteness $v_k(G) \leq m,$
 where $1\leq m \leq n-k.$ Then the following holds.
  $$RDD(G)\leq\frac{2n^3-n^2-3mn+2m}{2}-\frac{(n-m)(3n-s-1)s}{2}-\frac{t(s+1)(3n-2s-2)}{2}$$
with equality holds if and only if $G\cong K_m\vee((k-t)\overline{K_s}\vee t\overline{K_{s+1}}),$\\
where $(k-t)\overline{K_s}=\underbrace{\overline{K_s}\vee\ldots\vee\overline{K_s}}_{k-t},$and
$t\overline{K_{s+1}}=\underbrace{\overline{K_{s+1}}\vee\ldots\vee\overline{K_{s+1}}}_t.$
\end{theorem}

{\bf Proof.} Since the reciprocal degree distance  is a monotonic
increasing topological index, by Theorem $3.2,$
there exist $k$ non-negative integers $s_1,s_2,\ldots,s_k$
satisfying $s_1+s_2+\ldots+s_k=n-m,$
such that $\widehat{G}=K_m\vee\big(\overline{K_{s_1}}\vee\overline{K_{s_2}}\vee\ldots\vee\overline{K_{s_k}}\big)$
is the graph in $\mathscr {G}_{n,m,k}$ with the maximum reciprocal degree distance.

In the following, we will determine the values of $s_1,s_2,\ldots,s_k.$

$RDD(K_m\vee(\overline{K_{s_1}} \vee\overline{K_{s_2}} \vee\ldots\vee \overline{K_{s_k}}))$\\
=$\sum\limits_{\{u,v\}\subseteq V(\widehat{G})}\frac{d(u)+d(v)}{d(u,v)}$\\
=$C_m^2\times \frac{(n-1)+(n-1)}{1}+\sum\limits_{i=1}^{k} C_{s_i}^2\times \frac{(n-s_i)+(n-s_i)}{2} \\
+\sum\limits_{i=1}^{k}m\times s_i\times\frac{(n-s_i)+(n-1)}{1}
+\sum\limits_{1\leq i<j\leq k}s_is_j\times \frac{(n-s_i)+(n-s_j)}{1}.$

We claim that $\widehat{G}$ is the graph in $\mathscr {G}_{n,m,k}$ with the maximum reciprocal degree distance
when $s_i = s_j $ or $|s_i-s_j|=1,$ for all $1\leq i,j\leq k.$
If there exist $1\leq i,j\leq k,$ such that $|s_i-s_j|\geq 2.$
Wihtout loss of generality, we assume that $s_1-s_2\geq2.$
By moving one vertex from the part of $\overline{K_{s_1}}$ to the part of  $\overline{K_{s_2}},$
we get a new graph
$\widetilde{G}=K_m\vee(\overline{K_{s_1-1}} \vee\overline{K_{s_2+1}}\vee\ldots\vee\overline{K_{s_k}}),$
which is also in $\mathscr{G}_{n,m,k}.$

Then we obtain that\\
$RDD(\widehat{G})-RDD(\widetilde{G})$\\
=$\frac{s_1(s_1-1)(n-s_1)}{2}+\frac{s_2(s_2-1)(n-s_2)}{2}+ms_1(2n-s_1-1)+ms_2(2n-s_2-1)\\
+\sum\limits_{i=2}^{k}s_1s_i(2n-s_1-s_i)+\sum\limits_{i=3}^{k}s_2s_i(2n-s_2-s_i)\\
-\frac{(s_1-1)(s_1-2)(n-s_1+1)}{2}-\frac{(s_2+1)s_2(n-s_2+1)}{2}\\
-m(s_1-1)(2n-s_1)-m(s_2+1)(2n-s_2-2)\\
-\sum\limits_{i=3}^{k}(s_1-1)s_i(2n-s_1-s_i+1)-(s_1-1)(s_2+1)(2n-s_1-s_2)\\
-\sum\limits_{i=3}^{k}(s_2+1)s_i(2n-s_2-s_i-1)$\\
=$-(6n-2-3s_1-3s_2)(s_1-s_2-1).$

Since $s_1-s_2\geq2,$ then $RDD(\widehat{G})< RDD(\widetilde{G}),$
which is a contradiction.

Therefore, for any $1\leq i,j\leq k,$ then $|s_i-s_j|\leq 1.$ For $s_1+s_2+\ldots +s_k=n-m=sk+t,$
we have $\overline{K_{s_1}} \vee\overline{K_{s_2}} \vee\ldots\vee\overline{K_{s_k}}=(k-t)\overline{K_s}\vee t\overline{K_{s+1}},$\\
where $(k-t)\overline{K_s}=\underbrace{\overline{K_s}\vee\ldots\vee\overline{K_s}}_{k-t},$ and
$t\overline{K_{s+1}}=\underbrace{\overline{K_{s+1}}\vee\ldots\vee\overline{K_{s+1}}}_t.$

Then we obtain that
$RDD(\widehat{G})=\frac{2n^3-n^2-3mn+2m}{2}-\frac{(n-m)(3n-s-1)s}{2}-\frac{t(s+1)(3n-2s-2)}{2}.$\qed

\begin{coro}
Let $G\in\mathscr{G}_{n,m,2}.$ Then the following holds.

(1) If $n-m$ is even, then
$$RDD(G) \leq  \frac{3n^3-m^3-3m^2n+9mn^2-2n^2+2m^2-16mn+8m}{8},$$
with equality holds if and only if
 $$G\cong K_{m}\vee\big(\overline{K_{\frac{n-m}{2}}}\vee\overline{K_{\frac{n-m}{2}}}\big).$$

(2) If $n-m$ is odd, then
 $$ RDD(G) \leq \frac{3n^3-m^3-3m^2n+9mn^2-2n^2+2m^2-16mn-3n+5m+2}{8},$$
with equality holds if and only if
$$G\cong K_{m}\vee\big(\overline{K_{\frac{n-m+1}{2}}}\vee\overline{K_{\frac{n-m-1}{2}}}\big).$$

\end{coro}

\subsection{The minimum eccentricity distance sum  in $\mathscr{G}_{n,m,k}$}
\begin{theorem}
Let $G$ be a connected simple graph of order $n$ with vertex $k$-partiteness $v_k(G) \leq m,$
 where $1\leq m \leq n-k.$ Then the following holds.
  $$\xi^d(G)\geq m(n-1)+2t(s+1)(n+s-1)+2(k-t)s(n+s-2)$$
with equality holds if and only if $G\cong K_m\vee((k-t)\overline{K_s}\vee t\overline{K_{s+1}}),$\\
where $(k-t)\overline{K_s}=\underbrace{\overline{K_s}\vee\ldots\vee\overline{K_s}}_{k-t},$ and
$t\overline{K_{s+1}}=\underbrace{\overline{K_{s+1}}\vee\ldots\vee\overline{K_{s+1}}}_t.$
\end{theorem}

{\bf Proof.} Since the eccentricity distance sum  is a monotonic
decreasing topological index, by Theorem $3.1,$
there exist $k$ non-negative integers $s_1,s_2,\ldots,s_k$
satisfying $s_1+s_2+\ldots+s_k=n-m,$
such that $\widehat{G}=K_m\vee\big(\overline{K_{s_1}}\vee\overline{K_{s_2}}\vee\ldots\vee\overline{K_{s_k}}\big)$
is the graph in $\mathscr {G}_{n,m,k}$ with the minimum eccentricity distance sum.

In the following, we will determine the values of $s_1,s_2,\ldots,s_k.$

$\xi^d(K_m\vee(\overline{K_{s_1}} \vee\overline{K_{s_2}} \vee\ldots\vee \overline{K_{s_k}}))$\\
=$\sum\limits_{u\in V(\widehat{G})}\varepsilon(u)D(u)$\\
=$m(n-1)+\sum\limits_{i=1}^{k}2{s_i}^2\times (n+s_i-2).$

We claim that $\widehat{G}$ is the graph in $\mathscr {G}_{n,m,k}$ with the minimum eccentricity distance sum
when $s_i = s_j $ or $|s_i-s_j|=1,$ for all $1\leq i,j\leq k.$
If there exist $1\leq i,j\leq k$, such that $|s_i-s_j|\geq 2.$
Without loss of generality, we assume that $s_1-s_2\geq2.$
By moving one vertex from the part of $\overline{K_{s_1}}$ to the part of  $\overline{K_{s_2}},$
we get a new graph
$\widetilde{G}=K_m\vee(\overline{K_{s_1-1}} \vee\overline{K_{s_2+1}}\vee\ldots\vee\overline{K_{s_k}}),$
which is also in $\mathscr{G}_{n,m,k}.$

Then we obtain that\\
$\xi^d(\widehat{G})-\xi^d(\widetilde{G})$\\
$=2s_1(n+s_1-2)+2s_2(n+s_2-2)-2(s_1-1)(n+s_1-3)-2(s_2+1)(n+s_2-1)$\\
$=4(s_1-s_2-1).$

Since $s_1-s_2\geq2,$ then $\xi ^d(\widehat{G})> \xi^d(\widetilde{G}),$
which is a contradiction.

Therefore, for any $1\leq i,j\leq k,$ then $ |s_i-s_j|\leq 1.$ For $s_1+s_2+\ldots +s_k=n-m=sk+t,$
we have  $\overline{K_{s_1}} \vee\overline{K_{s_2}} \vee\ldots\vee\overline{K_{s_k}}=(k-t)\overline{K_s}\vee t\overline{K_{s+1}},$\\
where $(k-t)\overline{K_s}=\underbrace{\overline{K_s}\vee\ldots\vee\overline{K_s}}_{k-t},$ and
$t\overline{K_{s+1}}=\underbrace{\overline{K_{s+1}}\vee\ldots\vee\overline{K_{s+1}}}_t.$

Then we obtain that\\
$\xi^d(\widehat{G})
=m(n-1)+2t(s+1)(n+s-1)+2(k-t)s(n+s-2).$\qed

\begin{coro}
Let $G\in\mathscr{G}_{n,m,2}.$ Then the following holds.

(1) If $n-m$ is even, then
$$\xi^d(G) \geq 3n^2-4n-3mn+m^2+3m,$$
with equality holds if and only if
 $$G\cong K_{m}\vee\big(\overline{K_{\frac{n-m}{2}}}\vee\overline{K_{\frac{n-m}{2}}}\big).$$

(2) If $n-m$ is odd, then
 $$ \xi^d(G) \geq  3n^2-4n-3mn+m^2+3m-1 ,$$
with equality holds if and only if
$$G\cong K_{m}\vee\big(\overline{K_{\frac{n-m+1}{2}}}\vee\overline{K_{\frac{n-m-1}{2}}}\big).$$

\end{coro}

\subsection{The  minimum adjacent eccentric distance sum index  in $\mathscr{G}_{n,m,k}$}
\begin{theorem}
Let $G$ be a connected simple graph of order $n$ with vertex $k$-partiteness $v_k(G) \leq m,$
 where $1\leq m \leq n-k.$ Then the following holds.
  $$\xi^{ad}(G)\geq m+\frac{2t(s+1)(n+s-1)}{n-s-1}+\frac{2(k-t)s(n+s-2)}{n-s}$$
with equality holds if and only if $G\cong K_m\vee((k-t)\overline{K_s}\vee t\overline{K_{s+1}}),$\\
where $(k-t)\overline{K_s}=\underbrace{\overline{K_s}\vee\ldots\vee\overline{K_s}}_{k-t},$ and
$t\overline{K_{s+1}}=\underbrace{\overline{K_{s+1}}\vee\ldots\vee\overline{K_{s+1}}}_t.$
\end{theorem}

{\bf Proof.} Since the adjacent eccentric distance sum index  is a monotonic
decreasing topological index, by Theorem $3.1,$
there exist $k$ non-negative integers $s_1,s_2,\ldots,s_k$
satisfying $s_1+s_2+\ldots+s_k=n-m,$
such that $\widehat{G}=K_m\vee\big(\overline{K_{s_1}}\vee\overline{K_{s_2}}\vee\ldots\vee\overline{K_{s_k}}\big)$
is the graph in $\mathscr {G}_{n,m,k}$ with the minimum adjacent eccentric distance sum index.

In the following, we will determine the values of $s_1,s_2,\ldots,s_k.$

$\xi^{ad}(K_m\vee(\overline{K_{s_1}} \vee\overline{K_{s_2}} \vee\ldots\vee \overline{K_{s_k}}))$\\
=$\sum\limits_{u\in V(\widehat{G})}\frac{\varepsilon(u)D(u)}{d(u)}$\\
$=m\frac{n-1}{n-1}+\sum\limits_{i=1}^{k}s_i\frac{2(n+s_i-2)}{n-s_i}.$\\

We claim that $\widehat{G}$ is the graph in $\mathscr {G}_{n,m,k}$ with the minimum adjacent eccentric distance sum index
when $s_i = s_j $ or $|s_i-s_j|=1,$ for all $1\leq i,j\leq k.$
If there exist $1\leq i,j\leq k,$ such that $|s_i-s_j|\geq 2.$
Without loss of generality, we assume that $s_1-s_2\geq2.$
By moving one vertex from the part of $\overline{K_{s_1}}$ to the part of  $\overline{K_{s_2}},$
we get a new graph
$\widetilde{G}=K_m\vee(\overline{K_{s_1-1}} \vee\overline{K_{s_2+1}}\vee\ldots\vee\overline{K_{s_k}}),$
which is also in $\mathscr{G}_{n,m,k}.$

Then we obtain that\\
$\xi^{ad}(\widehat{G})-\xi^{ad}(\widetilde{G})$\\
$=\frac{2s_1(n+s_1-2)}{n-s_1}+\frac{2s_2(n+s_2-2)}{n-s_2}-\frac{2(s_1-1)(n+s_1-3)}{n-s_1+1}-\frac{2(s_2+1)(n+s_2-1)}{n-s_2-1}$\\
$=\frac{4n(s_1-s_2-1)(2n^2+2n(s_1-s_2)+s_1+s_2)+6n^2(s_2+1)+2n^3}{(n-s_1)(nis_2)(n-s_1+1)(n-s_2-1)}.$

Since $s_1-s_2\geq2,$ then $\xi ^{ad}(\widehat{G})> \xi^{ad}(\widetilde{G}),$
which is a contradiction.

Therefore, for any $1\leq i,j\leq k,$ then $ |s_i-s_j|\leq 1.$ For $s_1+s_2+\ldots +s_k=n-m=sk+t,$
we have $\overline{K_{s_1}} \vee\overline{K_{s_2}} \vee\ldots\vee\overline{K_{s_k}}=(k-t)\overline{K_s}\vee t\overline{K_{s+1}},$\\
where $(k-t)\overline{K_s}=\underbrace{\overline{K_s}\vee\ldots\vee\overline{K_s}}_{k-t},$ and
$t\overline{K_{s+1}}=\underbrace{\overline{K_{s+1}}\vee\ldots\vee\overline{K_{s+1}}}_t.$

Then we obtain that\\
$\xi^{ad}(\widehat{G})
=m+\frac{2t(s+1)(n+s-1)}{n-s-1}+\frac{2(k-t)s(n+s-2)}{n-s}.$\qed

\begin{coro}
Let $G\in\mathscr{G}_{n,m,2}.$ Then the following holds.

(1) If $n-m$ is even, then
$$\xi^{ad}(G) \geq m+\frac{2(n-m)(3n-m-4)}{n+m},$$
with equality holds if and only if
 $$G\cong K_{m}\vee\big(\overline{K_{\frac{n-m}{2}}}\vee\overline{K_{\frac{n-m}{2}}}\big).$$

(2) If $n-m$ is odd, then
 $$ \xi^{ad}(G) \geq  m+\frac{(n-m+1)(3n-m-3)}{n+m-1}+\frac{(n-m-11)(3n-m-5)}{n+m+1} ,$$
with equality holds if and only if
$$G\cong K_{m}\vee\big(\overline{K_{\frac{n-m+1}{2}}}\vee\overline{K_{\frac{n-m-1}{2}}}\big).$$

\end{coro}

\subsection{The  maximum connective eccentricity index  in $\mathscr{G}_{n,m,k}$}
\begin{theorem}
Let $G$ be a connected simple graph of order $n$ with vertex $k$-partiteness $v_k(G) \leq m,$
 where $1\leq m \leq n-k.$ Then the following holds.
  $$\xi^{ce}(G)\geq m(n-1)+\frac{t(s+1)(n-s-1)}{2}+\frac{(k-t)s(n-s)}{2}$$
with equality holds if and only if $G\cong K_m\vee((k-t)\overline{K_s}\vee t\overline{K_{s+1}}),$\\
where $(k-t)\overline{K_s}=\underbrace{\overline{K_s}\vee\ldots\vee\overline{K_s}}_{k-t},$ and
$t\overline{K_{s+1}}=\underbrace{\overline{K_{s+1}}\vee\ldots\vee\overline{K_{s+1}}}_t.$
\end{theorem}

{\bf Proof.} Since the connective eccentricity index  is a monotonic
increasing topological index, by Theorem $3.2,$
there exist $k$ non-negative integers $s_1,s_2,\ldots,s_k$
satisfying $s_1+s_2+\ldots+s_k=n-m,$
such that $\widehat{G}=K_m\vee\big(\overline{K_{s_1}}\vee\overline{K_{s_2}}\vee\ldots\vee\overline{K_{s_k}}\big)$
is the graph in $\mathscr {G}_{n,m,k}$ with the maximum connective eccentricity index.

In the following, we will determine the values of $s_1,s_2,\ldots,s_k.$

$\xi^{ce}(K_m\vee(\overline{K_{s_1}} \vee\overline{K_{s_2}} \vee\ldots\vee \overline{K_{s_k}}))$\\
=$\sum\limits_{u\in V(\widehat{G})}\frac{d(u)}{\varepsilon (u)}$\\
=$m\frac{n-1}{1}+\sum\limits_{i=1}^{k}s_i\frac{(n-s_i)}{2}.$

We claim that $\widehat{G}$ is the graph in $\mathscr {G}_{n,m,k}$ with the maximum connective eccentricity index
when $s_i = s_j $ or $|s_i-s_j|=1,$ for all $1\leq i,j\leq k.$
If there exist $1\leq i,j\leq k,$ such that $|s_i-s_j|\geq 2.$
Without loss of generality, we assume that $s_1-s_2\geq2.$
By moving one vertex from the part of $\overline{K_{s_1}}$ to the part of  $\overline{K_{s_2}}$
we get a new graph
$\widetilde{G}=K_m\vee(\overline{K_{s_1-1}} \vee\overline{K_{s_2+1}}\vee\ldots\vee\overline{K_{s_k}}),$
which is also in $\mathscr{G}_{n,m,k}.$

Then we obtain that\\
$\xi^{ce}(\widehat{G})-\xi^{ce}(\widetilde{G})$\\
$=\frac{s_1(n-s_1)}{2}+\frac{s_2(n-s_2)}{2}-\frac{(s_1-1)(n-s_1+1)}{2}-\frac{(s_2+1)(n-s_2-1)}{2}$\\
$=-(s_1-s_2-1).$

Since $s_1-s_2\geq2,$ then $\xi ^{ce}(\widehat{G})< \xi^{ce}(\widetilde{G}),$
which is a contradiction.

Therefore, for any $1\leq i,j\leq k,$ then $ |s_i-s_j|\leq 1.$ For $s_1+s_2+\ldots +s_k=n-m=sk+t,$
we have $\overline{K_{s_1}} \vee\overline{K_{s_2}} \vee\ldots\vee\overline{K_{s_k}}=(k-t)\overline{K_s}\vee t\overline{K_{s+1}},$\\
where $(k-t)\overline{K_s}=\underbrace{\overline{K_s}\vee\ldots\vee\overline{K_s}}_{k-t},$and
$t\overline{K_{s+1}}=\underbrace{\overline{K_{s+1}}\vee\ldots\vee\overline{K_{s+1}}}_t.$

Then we obtain that\\
$\xi^{ce}(\widehat{G})
=m(n-1)+\frac{t(s+1)(n-s-1)}{2}+\frac{(k-t)s(n-s)}{2}.$\qed

\begin{coro}
Let $G\in\mathscr{G}_{n,m,2}.$ Then the following holds.

(1) If $n-m$ is even, then
$$\xi^{ce}(G) \geq m(n-1)+\frac{(n^2-m^2)}{4},$$
with equality holds if and only if
 $$G\cong K_{m}\vee\big(\overline{K_{\frac{n-m}{2}}}\vee\overline{K_{\frac{n-m}{2}}}\big).$$

(2) If $n-m$ is odd, then
 $$ \xi^{ce}(G) \geq  m(n-1)+\frac{(n^2-m^2-1)}{4} ,$$
with equality holds if and only if
$$G\cong K_{m}\vee\big(\overline{K_{\frac{n-m+1}{2}}}\vee\overline{K_{\frac{n-m-1}{2}}}\big).$$

\end{coro}

\subsection{The  maximum Zagreb indices  in $\mathscr{G}_{n,m,k}$}
\begin{theorem}
Let $G$ be a connected simple graph of order $n$ with vertex $k-$partiteness $v_k(G) \leq m,$
 where $1\leq m \leq n-k.$ Then the following holds.
  $$ M_1(G)\leq m(n-1)^2+t(s+1)(n-s-1)^2+(k-t)s(n-s)^2$$
  $$ M_2(G)\leq \frac {m(m-1)(n-1)^2}{2}+m(n-1)t(s+1)(n-s-1)+m(n-1)(k-t)s(n-s)$$
  $$ \Pi_1(G)\leq (n-1)^{2m}(n-s)^{2s(k-t)}(n-s-1)^{2t(s+1)}$$
  $$ \Pi_2(G)\leq (n-1)^{m(n-1)}(n-s)^{s(n-s)(k-t)}(n-s-1)^{(s+1)(n-s-1)t}$$
with equality holds if and only if $G\cong K_m\vee((k-t)\overline{K_s}\vee t\overline{K_{s+1}}),$\\
where $(k-t)\overline{K_s}=\underbrace{\overline{K_s}\vee\ldots\vee\overline{K_s}}_{k-t},$and
$t\overline{K_{s+1}}=\underbrace{\overline{K_{s+1}}\vee\ldots\vee\overline{K_{s+1}}}_t.$
\end{theorem}

{\bf Proof.} Since the Zagreb indices  are monotonic
increasing topological indices, by Theorem $3.2,$
there exist $k$ non-negative integers $s_1,s_2,\ldots,s_k$
satisfying $s_1+s_2+\ldots+s_k=n-m,$
such that $\widehat{G}=K_m\vee\big(\overline{K_{s_1}}\vee\overline{K_{s_2}}\vee\ldots\vee\overline{K_{s_k}}\big)$
is the graph in $\mathscr {G}_{n,m,k}$ with the maximum Zagreb indices.

In the following, we will determine the values of $s_1,s_2,\ldots,s_k.$\\
$M_1(K_m\vee(\overline{K_{s_1}} \vee\overline{K_{s_2}} \vee\ldots\vee \overline{K_{s_k}}))$\\
=$\sum\limits_{u\in V(\widehat{G})}d(u)^2$\\
=$m(n-1)^2+\sum\limits_{i=1}^{k}s_i(n-s_i)^2.$\\
$M_2(K_m\vee(\overline{K_{s_1}} \vee\overline{K_{s_2}} \vee\ldots\vee \overline{K_{s_k}}))$\\
=$\sum\limits_{uv\in E(\widehat{G})}d(u)d(v)$\\
=$\frac{m(m-1)(n-1)^2}{2}+\sum\limits_{i=1}^{k}ms_i(n-1)(n-s_2)+\sum\limits_{1\leq i <j\leq k}s_is_j(n-s_i)(n-s_j).$\\
$\Pi_1(K_m\vee(\overline{K_{s_1}} \vee\overline{K_{s_2}} \vee\ldots\vee \overline{K_{s_k}}))$\\
=$\prod \limits_{u\in V(\widehat{G})}d(u)^2$\\
=$(n-1)^{2m}\prod \limits_{i=1}^{k}(n-s_i)^{2s_i}.$\\
$\Pi_2(K_m\vee(\overline{K_{s_1}} \vee\overline{K_{s_2}} \vee\ldots\vee \overline{K_{s_k}}))$\\
=$\prod \limits_{u\in V(\widehat{G})}d(u)^{d(u)}$\\
=$(n-1)^{m(n-1)}\prod \limits_{i=1}^{k}(n-s_i)^{s_i(n-s_i)}.$

We claim that $\widehat{G}$ is the graph in $\mathscr {G}_{n,m,k}$ with the maximum Zagreb indices
when $s_i = s_j $ or $|s_i-s_j|=1,$ for all $1\leq i,j\leq k.$
If there exist $1\leq i,j\leq k,$ such that $|s_i-s_j|\geq 2.$
Without loss of generality, we assume that $s_1-s_2\geq2.$
By moving one vertex from the part of $\overline{K_{s_1}}$ to the part of  $\overline{K_{s_2}},$
we get a new graph
$\widetilde{G}=K_m\vee(\overline{K_{s_1-1}} \vee\overline{K_{s_2+1}}\vee\ldots\vee\overline{K_{s_k}}),$
which is also in $\mathscr{G}_{n,m,k}.$

Then we obtain that\\
$M_1(\widehat{G})-M_1(\widetilde{G})$\\
$=s_1(n-s_1)^2+s_2(n-s_2)^2-(s_1-1)(n-s_1+1)^2-(s_2+1)(n-s_2-1)^2$\\
$=-(s_1-s_2-1)(4n-3s_1-3s_2).$\\
$M_2(\widehat{G})-M_2(\widetilde{G})$\\
$=ms_1(n-1)(n-s_1)+ms_2(n-1)(n-s_2)+\sum\limits_{i=2}^{k}s_1s_i(n-s_1)(n-s_i)+\sum\limits_{i=3}^{k}s_2s_i(n-s_2)(n-s_i)$\\
$-m(s_1-1)(n-1)(n-s_1+1)-m(s_2+1)(n-1)(n-s_2-1)-(s_1-1)(s_2+1)(n-s_1+1)(n-s_2-1)$\\
$-\sum\limits_{i=3}^{k}(s_1-1)s_i(n-s_1+1)(n-s_i)-\sum\limits_{i=3}^{k}(s_2+1)s_i(n-s_2-1)(n-s_i)$\\
$=-(s_1-s_2-1)(2m(n-1)+s_1s_2+(n-s_1)(n-s_2)+s_1-s_2-1+\sum\limits_{i=3}^{k}2s_i(n-s_i)).$\\
$\frac{\Pi_1(\widehat{G})}{\Pi_1(\widetilde{G})}$\\
$=\frac{(n-s_1)^{2s_1}(n-s_2)^{2s_2}}{(n-s_1+1)^{2(s_1-1)}(n-s_2-1)^{2(s_2+1)}}$\\
$=(\frac{n-s_1}{n-s_2-1})^{2s_1}(\frac{n-s_2-1}{n-s_1+1})^{2(s_1-s_2-1)}(\frac{n-s_2}{n-s_1+1})^{2s_2}.$\\
$\frac{\Pi_2(\widehat{G})}{\Pi_2(\widetilde{G})}$\\
$=\frac{(n-s_1)^{s_1(n-s_1)}(n-s_2)^{s_2(n-s_2)}}{(n-s_1+1)^{(s_1-1)(n-s_1+1}(n-s_2-1)^{(s_2+1)(n-s_2-1)}}$\\
$=(\frac{n-s_1}{n-s_1+1})^{s_1(n-s_1)}(\frac{n-s_2}{n-s_2-1})^{s_2(n-s_2)}\frac{(n-s_1+1)^(n-2s_1+1)}{(n-s_2-1)^(n-2s_2+1)}.$

Since $s_1-s_2\geq2,$ then $M_1(\widehat{G})< M_1(\widetilde{G}),$ $M_2(\widehat{G})< M_2(\widetilde{G}),$ $\Pi_1(\widehat{G})< \Pi_1(\widetilde{G}),$ $\Pi_2(\widehat{G})< \Pi_2(\widetilde{G}).$
which is a contradiction.

Therefore, for any $1\leq i,j\leq k,$ then $ |s_i-s_j|\leq 1.$ For $s_1+s_2+\ldots +s_k=n-m=sk+t,$
we have  $\overline{K_{s_1}} \vee\overline{K_{s_2}} \vee\ldots\vee\overline{K_{s_k}}=(k-t)\overline{K_s}\vee t\overline{K_{s+1}},$\\
where $(k-t)\overline{K_s}=\underbrace{\overline{K_s}\vee\ldots\vee\overline{K_s}}_{k-t},$ and
$t\overline{K_{s+1}}=\underbrace{\overline{K_{s+1}}\vee\ldots\vee\overline{K_{s+1}}}_t.$

Then we obtain that\\
$M_1(\widehat{G})
=m(n-1)^2+t(s+1)(n-s-1)^2+(k-t)s(n-s)^2.$\\
$M_2(\widehat{G})
=\frac{m(m-1)(n-1)^2}{2}+m(n-1)t(s+1)(n-s-1)+m(n-1)(k-t)s(n-s).$\\
$\Pi_1(\widehat{G})
=(n-1)^{2m}(n-s)^{2s(k-t)}(n-s-1)^{2(s+1)t}.$\\
$\Pi_2(\widehat{G})
=(n-1)^{m(n-1)}(n-s)^{s(n-s)(k-t)}(n-s-1)^{(s+1)(n-s-1)t}.$
\qed

\begin{coro}
Let $G\in\mathscr{G}_{n,m,2}.$ Then the following holds.

(1) If $n-m$ is even, then
$$ M_1(G) \leq m(n-1)^2+\frac{(n^3+n^2m-m^2n-m^3)}{4},$$
$$ M_2(G) \leq \frac{m(m-1)(n-1)^2}{2}+\frac{(n^2-m^2)(n^2+8mn-m^2-8m)}{16},$$
$$ \Pi_1(G) \leq (n-1)^{2m}(\frac{n+m}{2})^{2(n-m)} ,$$
$$ \Pi_2(G) \leq (n-1)^{m(n-1)}(\frac{n+m}{2})^{\frac{n^2-m^2}{2}} ,$$
with equality holds if and only if
$$G\cong K_{m}\vee\big(\overline{K_{\frac{n-m}{2}}}\vee\overline{K_{\frac{n-m}{2}}}\big).$$

(2) If $n-m$ is odd, then
$$ M_1(G) \leq  m(n-1)^2+\frac{(n^3+n^2m-m^2n-n-m^3-3m)}{4} ,$$
$$ M_2(G) \leq \frac{m(m-1)(n-1)^2}{2}+\frac{[((n+m)^2-1))((n-m)^2-1)+4m(n-1)(2n^2-2m^2-2))}{16},$$
$$ \Pi_1(G) \leq (n-1)^{2m}(\frac{n+m+1}{2})^{n-m-1}(\frac{n+m-1}{2})^{n-m+1} ,$$
$$ \Pi_2(G) \leq (n-1)^{m(n-1)}(\frac{n+m-1}{2})^{\frac{n^2-m^2+2m-1}{4}}(\frac{n+m+1}{2})^{\frac{n^2-m^2-2m-1}{4}} ,$$
with equality holds if and only if
$$G\cong K_{m}\vee\big(\overline{K_{\frac{n-m+1}{2}}}\vee\overline{K_{\frac{n-m-1}{2}}}\big).$$

\end{coro}


\begin{thebibliography}{99}
\bibitem{ash}  A.R. Ashrafi, M. Saheli, M. Ghorbani, The eccentric connectivity index of nanotubes and nanotori, {\it J. Comput. Appl. Math.} 235 (2011) 4561-4566.
\vspace{-0.3cm}
\bibitem{Bondy} J.A. Bondy,U.S.R. Murty, Graph theory with applications, {\it Macmillan, London and Elsevier, New York.} 1976.
\vspace{-0.3cm}
\bibitem{Fallat} S. Fallat, Y.Z. Fan, Bipartiteness and the least eigenvalue of signless Laplacian of graphs, {\it Linear Algebra Appl.} 436 (2012) 3254-3267.
\vspace{-0.3cm}
\bibitem{Freitas2016} M. A. A. de Freitas, I. Gutman, M. Robbiano, Graphs with maximum Laplacian-energy-like invariant and incidence
energy, {\it MATCH Commun. Math. Comput. Chem.} 75 (2016) 331-342.
\vspace{-0.3cm}
\bibitem{gup2000} S. Gupta, M. Singh, A.K. Madan, Connective eccentricity index:a novel topological descriptor for predicting biological activity, {\it J. Mol. Graph. Model.} 18 (2000) 18-25.
\vspace{-0.3cm}
\bibitem{gup2002} S. Gupta, M. Singh, A.K. Madan, Application of graph theory:relationship of eccentric connectivity index and Wiener's index with anti-inflammatory activity, {\it J. Math. Anal. Appl.} 266 (2002) 259-268.
\vspace{-0.3cm}
\bibitem{gut1975} I. Gutman, B. Ru\u{s}\u{c}i\'{c}, N. Trinajsti\'{c}, C.F. Wilcox, Graph theory and molecular orbitals, {\it J. Chen. Phys.} 62 (1975) 3390-3405.
\vspace{-0.3cm}
\bibitem{gut1972} I. Gutman, N. Trinajsti\'{c}, Graph theory and molecular orbitals.Total $\pi$-electron energy of alternant hydrocarbons, {\it Chen. Phys. Lett.} 17 (1972) 535-538.
\vspace{-0.3cm}
\bibitem{gut2013} I. Gutman, Degree-based topological indices, {\it Croatica Chemica.} 86 (2013) 351-361.
\vspace{-0.3cm}
\bibitem{gut2016} I. Gutman, L. Medina C, P. Pizarro, M. Robbiano, Graphs with maximum Laplacian and signless Laplacian Estrada index, {\it Discrete Appl. Math.} 339 (2016) 2664-2671.
\vspace{-0.3cm}
\bibitem{hong2017} H.S. Li, S.C. Li, H.H. Zhang, On the maximal connective eccentricity index of bipartite graphs with some given parameters, {\it J. Math. Anal. Appl.}  454 (2017) 453-467.
\vspace{-0.3cm}
\bibitem{hua2}H.B. Hua, S.G. Zhang, On the reciprocal degree distance of graphs, {\it Discrete Appl. Math.} 160 (2012) 1152-1163.
\vspace{-0.3cm}
\bibitem{liu} J.B. Liu, X.F. Pan, Minimizing Kirchhoff index among graphs with a given vertex bipartiteness, {\it Appl. Math. Comput.} 291 (2016) 84-88.
\vspace{-0.3cm}
\bibitem{Ni} A. Nihat, K.C. Das, A.S. \c{C}evik, Some properties on the tensor product of graphs obtained by monogenic semigroups, {\it Appl. Math. Comput.} 235 (2014) 352-357.
\vspace{-0.3cm}
\bibitem{sar2002} S. Sardana, A.K. Madan, Predicting anti-HIV activity of TIBO derivatives: A computational approach using a novel topological descriptor, {\it Mol. Model} 8 (2002) 258-265.
\vspace{-0.3cm}
\bibitem{rob} M. Robbiano, K. Tapia Morales, B. San Mart\'{\i}n, Extremal graphs with bounded vertex bipartiteness number, {\it Linear Algebra Appl.} 493 (2016) 28-36.
\vspace{-0.3cm}
\bibitem{tod2010} R. Todeschini, V. Conaoni, New local vertex invariants and molecular descriptors based on funtions of the vertex degrees, {\it Math. Comput.} 64 (2010) 359-372.
\vspace{-0.3cm}
\bibitem{wie}  H. Wiener, Structural determination of paraffin boiling point, {\it J. Amer. Chem. Soc.} 69 (1947) 17-20.
\vspace{-0.3cm}

\end{thebibliography}
\end{document}